\documentclass{elsarticle}

\usepackage{amsmath}
\usepackage{amsfonts} 
\usepackage{bm}
\usepackage{tikz}
\usepackage{pgfplots}
\usepgflibrary{arrows}
\usepackage{rotating} 
\graphicspath{{pictures/}}

\newtheorem{thm}{Theorem}

\newproof{pf}{Proof}

\journal{arXiv} 

\begin{document}

\begin{frontmatter}

\title{Fundamental mode exact schemes for unsteady problems}

\author[nsi,rudn]{Petr N. Vabishchevich\corref{cor}}
\ead{vabishchevich@gmail.com}

\address[nsi]{Nuclear Safety Institute, Russian Academy of Sciences, Moscow, Russia}
\address[rudn]{Peoples' Friendship University of Russia (RUDN University), Moscow, Russia}

\cortext[cor]{Corresponding author}

\begin{abstract}
The problem of increasing the accuracy of an approximate solution is considered for boundary value problems 
for parabolic equations. For ordinary differential equations (ODEs),
nonstandard finite difference schemes are in common use for this problem. 
They are based on a modification of standard discretizations of time derivatives and, in some cases, allow
to obtain the exact solution of problems. For multidimensional problems, we can consider the problem of 
increasing the accuracy only for the most important components of the approximate solution. In the present work, 
new unconditionally stable schemes for parabolic problems are constructed, which are exact for the fundamental mode.
Such two-level schemes are designed via a modification of standard schemes with weights 
using Pad\'{e} approximations. Numerical results obtained for a model problem demonstrate advantages of
the proposed fundamental mode exact schemes.
\end{abstract}


\end{frontmatter}

\section{Introduction}

In numerical solving time-dependent boundary value problems, the emphasis is on using
computational algorithms of higher accuracy  (see, e.g., \cite{HundsdorferVerwer2003,Gustafsson2008}).
Along with increasing the accuracy of discretization in space,
increasing the accuracy of discretization in time is also considered 
taking into account numerical methods developed for ordinary differential
equations \cite{Ascher,LeVeque}. In view of features of unsteady
problems for partial differential equations, we, first of all, should be guided by 
methods of numerical solving Cauchy problems for stiff systems of ordinary differential equations
\cite{HairerWanner2010,Butcher2008}.

As a rule, discretization in time for numerical solving boundary value problems for parabolic equations 
is constructed on the basis of two- or three-level schemes. To study
stability of such difference schemes, the theory of stability and well-posedness
for operator-difference schemes \cite{Samarskii1989,SamarskiiMatusVabishchevich} is employed. 
In particular, stability conditions for standard schemes with weights ($\theta$-method) are obtained for
linear problems in different Hilbert spaces.
For two-level schemes, where the solution at two consecutive time levels is involved, 
polynomial approximations are  explicitly or implicitly employed for operators of difference schemes.
The Runge-Kutta methods (see \cite{Butcher2008,DekkerVerwer1984}) are well-known examples of such schemes
widely used in modern computational practice. The main feature of multilevel schemes (multistep methods) results in
the approximation of time derivatives with higher accuracy on a multipoint stencil. 
Multistep methods based on numerical backward differentiation formulas \cite{Gear1971}
should be mentioned as a typical example.

In some cases, it is possible to construct exact difference schemes.
The main possibilities of their construction for solving  Cauchy problems for ordinary differential 
equations are reflected in the works \cite{mickens1994nonstandard,mickens2002nonstandard,anguelov2001contributions}.
The idea of increasing accuracy by modifying approximations of time derivatives
can be realized not only for linear equations, but also for nonlinear ones
and for systems of equations. Separately, we highlight the works 
\cite{mickens1999nonstandard,anguelov2005non,patidar2005use,hernandez2013nonstandard}, 
where such a technology of constructing discretization in time is applied to solving parabolic equations.
Using perturbations with a low parameter for standard difference approximations, we can
provide a new quality of the solution, which is associated, for example, with the fulfillment of some conservation law.

For a Cauchy problem for a homogeneous linear parabolic equation, the solution can be represented
as a superposition of modes, which are associated with the corresponding eigenfunctions of the elliptic operator. 
Dominant modes are the most important, since higher modes decay very quickly.
Such a qualitative behavior of the solution must hold for a discrete problem as well.
For long times, the fundamental mode of the solution is highlighted, when a regular regime is considered (see, e.g., \cite{luikov2012analytical}).
In the theory of difference schemes \cite{samarskii1996computational}, there is recognized
a class of asymptotically stable difference schemes that provide the correct behavior of the approximate solution for large times.
SM (Spectral Mimetic) properties of operator-difference schemes for numerical solving
the Cauchy problem for evolutionary equations of first order are  associated with the
 time-evolution of individual modes of the solution. 
For long times, the fundamental mode of the solution is highlighted for a regular regime.
With this in mind, SM-stable difference schemes are constructed \cite{Vabischevich2010b}.

Here we construct difference schemes for parabolic equations, which are exact for the most important component 
of the solution, namely, the fundamental mode is calculated exactly. Two-level fundamental mode exact schemes 
(FMES) are constructed using Pad\'{e} approximations. 

The paper is organized as follows. 
A boundary value problem for a linear parabolic equation with a self-adjoint
elliptic operator of second order is considered in Section 2.
Using finite element approximations in space, we formulate
a Cauchy problem for the corresponding differential-operator equation.
In Section 3,  we consider standard two-level schemes with weights and formulate stability conditions.
Schemes that are exact for the solution associated with the fundamental eigenfunction of 
the corresponding elliptic operator are constructed in Section 4.
Section 5 presents the results of numerical experiments on studying FMES accuracy
for a model initial-boundary value parabolic problem.

\section{Problem statement}

Let $\Omega$ be a bounded domain ($\Omega \subset \mathbb{R}^d, d = 2,3$)
with a piecewise smooth boundary $\partial \Omega$.
Define an elliptic operator $\mathcal{A}$ such that
\begin{equation}\label{1}
 \mathcal{A} u :=  - \nabla (k(\bm x) \nabla u ) + c(\bm x) u ,
 \quad \bm x \in \Omega . 
\end{equation} 
The operator $A$ is defined on the set of functions $u({\bm x})$ that satisfy
on the boundary $\partial\Omega$ the following conditions:
\begin{equation}\label{2}
  k({\bm x}) \frac{\partial u }{\partial n } + \mu ({\bm x}) u = 0,
  \quad {\bm x} \in \partial \Omega .
\end{equation} 
Coefficients $k(\bm x)$, $c(\bm x)$ and $\mu ({\bm x})$ are smooth functions in $\overline{\Omega}$ and
\[
 k(\bm x) \geq \kappa > 0,
 \quad \mu(\bm x) \geq 0 ,
 \quad \bm x \in \Omega . 
\]
We consider the Cauchy problem
\begin{equation}\label{3}
 \frac{d u}{d t} + \mathcal{A} u = 0,
 \quad 0 < t \leq T ,  
\end{equation} 
\begin{equation}\label{4}
 u(0) = u^0,
\end{equation}
where, for example, $u^0(\bm x) \in L_2(\Omega)$ with the notation  $u(t) = u(\bm x,t)$. 

Let $(\cdot,\cdot), \|\cdot\|$ be the scalar product and norm in $H = L_2(\Omega)$, respectively:
\[
(u,v) := \int_\Omega u({\bm x})v({\bm x}) d {\bm x},
\quad \|u\| := (u,u)^{1/2}.
\]
Multiplying equation (\ref{3}) by $v(\bm x) \in H^1(\Omega)$ and integrating over the domain $\Omega$, 
we arrive at the equality 
\begin{equation}\label{5}
 \left( \frac{d u}{d t}, v \right ) + a(u,v) = 0, 
 \quad \forall  v \in H^1(\Omega) ,
 \quad 0 < t \leq T .
\end{equation} 
Here $a(\cdot, \cdot)$  is the following bilinear form:
\[
 a(u,v) := \int_{\Omega} ( k \nabla u \cdot \nabla v + c \, u \, v)  d \bm x 
 + \int_{\partial \Omega}  \mu  \, u \, v \, d \bm x,
\] 
where
\[
 a(u,v) = a(v,u),
 \quad a(u,u) \geq \delta \|u\|^2,
\] 
with a constant $\delta$.
In view of (\ref{4}), we put
\begin{equation}\label{6}
 (u(0),v) = (u^0, v), 
 \quad \forall  v \in H^1(\Omega) .
\end{equation} 
The variational (weak) formulation of the problem (\ref{1})--(\ref{4}) consists in finding
$u(\bm x, t) \in H^1(\Omega)$, $0 < t \leq T$ that satisfies (\ref{5}), (\ref{6}) with the condition (\ref{2}) 
on the boundary. 

For the solution of the problem (\ref{5}), (\ref{6}), we have the a priori estimate
\begin{equation}\label{7}
 \|u(t)\| \leq \exp(-\delta t) \|u^0\| .
\end{equation}
To show this, we put $v = u$ in (\ref{5}) and get
\[
 \|u\| \frac{d }{d t}\|u\| + a(u,u) = 0.
\]
Taking into account the lower bound for the bilinear form $a(\cdot, \cdot)$, we get
\[
 \frac{d }{d t} \|u(t)\| + \delta  \|u(t)\| \leq 0.
\]
From this inequality, in view of (\ref{6}), it follows that the estimate (\ref{7}) holds.

For numerical solving the initial-boundary value problem (\ref{3}), (\ref{4}), discretization in space 
is constructed using  the finite element method \cite{brenner,Thomee2006}. 
The weak formulation (\ref{5}), (\ref{6}) is employed.
Define the subspace of finite elements $V^h \subset H^1(\Omega)$ and the discrete elliptic operator $A$ as
\[
(A y, v) = a(y,v),
\quad \forall \ y,v \in V^h . 
\]
The operator $A$ acts on the finite dimensional space $V^h$ and
\begin{equation}\label{8}
A = A^* \geq \delta_h I ,
\end{equation} 
where $I$ is the identity operator.

For the problem (\ref{3}), (\ref{4}), we put into the correspondence the operator equation for $w(t) \in V^h$:
\begin{equation}\label{9}
 \frac{d w}{d t} + A w = 0, 
 \quad 0 < t \leq T, 
\end{equation} 
\begin{equation}\label{10}
 w(0) = w^0, 
\end{equation} 
where $w^0 = P u^0$ with $P$ denoting $L_2$-projection onto $V^h$.
Similarly to (\ref{7}), we establish the estimate 
\begin{equation}\label{11}
 \|w(t)\| \leq \exp(-\delta_h t) \|w^0\| 
\end{equation} 
for the solution of the problem (\ref{9}), (\ref{10}).

To solve the problem (\ref{9}), (\ref{10}), the separation of variables method is applied.
For the spectral problem
\[
A \varphi_k = \lambda_k \varphi_k
\] 
we have real eigenvalues
\[
\lambda_1 \leq \lambda_2 \leq ... \leq  \lambda_{M_h}.
\]
The eigenfunctions  $ \varphi_k, \ \|\varphi_k\| = 1, \ k = 1,2, ...,   M_h$ form a basis in $V^h$. 
Therefore, for $y \in V^h$, we have
\[
 y = \sum_{k=1}^{M_h} y_k  \varphi_k , 
 \quad y_k =  (y,\varphi_k),
 \quad  k = 1,2, ...,   M_h.
\]
In view of this, the solution of the problem (\ref{9}), (\ref{10}) is represented in the form
\begin{equation}\label{12}
 w(\bm x, t) = \sum_{k=1}^{M_h} w_k^0 \exp(- \lambda_k t) \varphi_k (\bm x). 
\end{equation}
With this representation, we get the a priori estimate (\ref{1}) with $\delta_h= \lambda_1$. 

Various effects arising from modes can be taken into account in constructing difference schemes.
Assume that the fundamental mode is separated from other modes, i.e., $\lambda_1 < \lambda_2$.
Then, for large times, a regular regime is established for the solution:
\[
 w(\bm x, t) \approx  w_1^0 \exp(- \lambda_1 t) \varphi_1 (\bm x). 
\] 
If such a property is preserved after discretization in time, then we speak of 
asymptotic stability of difference schemes \cite{samarskii1996computational}.
In the representation (\ref{12}), higher modes decay faster with time.
Such a monotone behavior of different modes is inherited in SM-stable difference schemes
\cite{Vabischevich2010b}.

It seems natural to construct discretizations in time that are more accurate
in reproducing the time dependence of dominant modes. In particular, we distinguish nonstandard difference schemes, 
which are exact for the fundamental mode of the solution, namely, fundamental mode exact schemes.

\section{Nonstandard schemes with weights} 

The standard approach to the solution of the problem (\ref{9}), (\ref{10})  is associated with using two-level schemes.
Let $\tau$ be a step of a uniform grid in time such that 
$y^n = y(t^n), \ t^n = n \tau$, $n = 0,1, ..., N, \ N\tau = T$.
For a constant weight parameter $\sigma \ (0 < \sigma \leq 1)$, we 
approximate equation (\ref{9}) by the following two-level scheme:
\begin{equation}\label{13}
 \frac{y^{n+1} - y^{n}}{\tau} + A (\sigma y^{n+1} + (1-\sigma) y^{n}) = 0 ,
  \quad n=0,1,... , N-1,
\end{equation} 
with the initial conditions
\begin{equation}\label{14}
  y^0 = w^0,
\end{equation} 
where $\sigma$ is a parameter (weight). 
If $\sigma = 0$, then  (\ref{13}), (\ref{14}) becomes the explicit scheme, for 
$\sigma = 1$, we obtain the fully implicit scheme, whereas  $\sigma = 0.5$ corresponds to the symmetric 
(Crank--Nicolson) scheme. 

Similarly to (\ref{12}), the solution of the problem  (\ref{13}), (\ref{14}) can be represented as
\begin{equation}\label{15}
 y^n = \sum_{k=1}^{M_h} y^n_k \varphi_k,
  \quad n=0,1,... , N . 
\end{equation} 
Substitution of  (\ref{15}) into (\ref{13}) results in the following relation for an individual mode:
\[
 y^{n+1}_k = r(\sigma, \lambda_k \tau) y^n_k,
 \quad k = 1, 2, ..., M_h,
 \quad n=0,1,... , N-1, 
\] 
where
\[
 r(\sigma, \eta) := \frac{1-(1-\sigma) \eta}{1+ \sigma \eta} .
\] 

The scheme (\ref{13}), (\ref{14}) belongs to the class of FMES, if the following condition holds:
\begin{equation}\label{16}
 y^{n+1}_1 = \exp(- \lambda_1 \tau) y^n_1 .
\end{equation}  
For any $\sigma$, the equation
\[
 r(\sigma, \eta_1) = \exp(- \eta_1),
 \quad \eta_1 = \lambda_1 \tau
\]
does not hold. Thus, for the schemes with weights (\ref{13}), (\ref{14}), there is not any
fundamental mode exact scheme.

\begin{figure}[tbp]
  \begin{center}
    \includegraphics[width=0.75\linewidth] {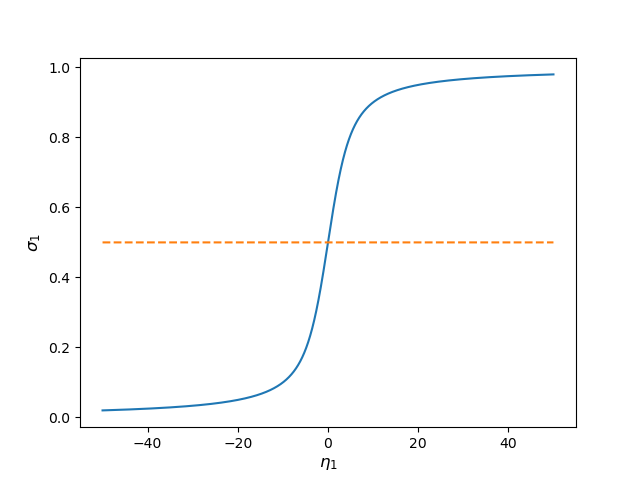}
	\caption{Weight value for FMES.}
	\label{fig:1}
  \end{center}
\end{figure} 

It is possible to select a weight value in such a way that the fundamental mode
is calculated exactly. The value $\sigma_1$ such that $r(\sigma_1, \eta_1) = \exp(- \eta_1)$
is shown in Fig.~\ref{fig:1}. We restrict ourselves to the values $|\eta_1| \ll 1$.
In the most important case, we have $\sigma_1 \approx 0.5$.
For computational practice, the fully implicit scheme ($\sigma = 1$) is more interesting,
but it cannot describe accurately the time-variation of the fundamental mode.

In our previous work \cite{afanas2013unconditionally}, special unconditionally stable schemes 
were constructed for convection-diffusion problems using a transition to new unknowns. 
The main feature of these schemes is connected with the negativity of the constant $\delta_h$ in the lower  bound
for the operator $A$ (the estimate (\ref{8}) in the present work).
Here we use this idea for constructing fundamental mode exact schemes.

For the solution of the problem (\ref{9}), (\ref{10}), we apply the representation
\begin{equation}\label{17}
 w(t) = \exp(- \lambda_1 t) v(t) .
\end{equation} 
Substitution of (\ref{17}) into (\ref{9}), (\ref{10}) results in the problem 
\begin{equation}\label{18}
 \frac{d v}{d t} + \widetilde{A} v = 0, 
 \quad 0 < t \leq T, 
\end{equation} 
\begin{equation}\label{19}
 v(0) = w^0, 
\end{equation} 
where
\begin{equation}\label{20}
 \widetilde{A} := A - \lambda_1 I , 
\end{equation} 
and
\[
 \widetilde{A} = \widetilde{A}^* \geq 0 .
\] 

For the problem (\ref{18}), (\ref{19}), we use the standard scheme with weights:
\begin{equation}\label{21}
 \frac{g^{n+1} - g^{n}}{\tau} + \widetilde{A} (\sigma g^{n+1} + (1-\sigma) g^{n}) = 0 ,
  \quad n=0,1,... , N-1,
\end{equation} 
\begin{equation}\label{22}
  g^0 = w^0 .
\end{equation} 
In view of (\ref{17}), suppose that $y^n = \exp(- \lambda_1 t^n) g^n$.
Thus, we arrive at the scheme
\begin{equation}\label{23}
 \frac{\exp(\lambda_1 \tau) y^{n+1} - y^{n}}{\tau} + \widetilde{A} (\sigma \exp(\lambda_1 \tau) y^{n+1} + (1-\sigma) y^{n}) = 0 ,
  \quad n=0,1,... , N-1,
\end{equation} 
with the initial condition
\begin{equation}\label{24}
  y^0 = w^0.
\end{equation}
This is a new scheme for the initial problem (\ref{9}), (\ref{10}).

A comparison of the nonstandard scheme (\ref{20}), (\ref{23}), (\ref{24})
with the standard scheme with weights (\ref{13}), (\ref{14}) indicates that
we modified not only discretization in time, but also the problem operator itself:
\begin{equation}\label{25}
 y^{n+1} \rightarrow \exp(\lambda_1 \tau) y^{n+1},
 \quad A \rightarrow A -  \lambda_1 I . 
\end{equation} 
Using the representation (\ref{15}), we obtain
\[
 \frac{\exp(\lambda_1 \tau) y_1^{n+1} - y_1^{n}}{\tau} = 0.
\] 
Thus, the condition (\ref{16}) is fulfilled and this scheme is of FMES-type.

The nonstandard scheme (\ref{20}), (\ref{23}), (\ref{24}) is unconditionally stable under 
the standard restriction on the weight $\sigma \geq 0.5$.
The proof is based on the well-known (see, e.g., \cite{Samarskii1989,SamarskiiMatusVabishchevich})
condition for the stability of the scheme with weights (\ref{21}), (\ref{22}).
Multiply equation (\ref{18}) scalarly by
\[
 \sigma g^{n+1} + (1-\sigma) g^{n} =
 \tau \left (\sigma - \frac{1}{2} \right ) \frac{g^{n+1} - g^{n}}{\tau} + \frac{1}{2} (g^{n+1} + g^{n}) .  
\] 
In view of $\widetilde{A} \geq 0$, we get
\[
 \tau \left (\sigma - \frac{1}{2} \right ) \left \| \frac{g^{n+1} - g^{n}}{\tau} \right \|^2 +
 \frac{1}{2 \tau } (\|g^{n+1}\|^2 - \|g^{n}\|^2) \leq 0 .
\] 
Thus, for $\sigma \geq 0.5$, we have
\[
 \|g^{n+1}\| \leq \|g^{n}\| \leq ... \leq \|w^{0}\| .
\]
Using the relation between $y^n$ and $ g^n$, we obtain for the scheme  (\ref{23}), (\ref{24}) the a priori estimate
\begin{equation}\label{26}
 \|y^{n+1}\| \leq \exp(- \lambda_1 t^{n+1}) \|w^{0}\| ,
 \quad n=0,1,... , N-1 . 
\end{equation}
This estimate is consistent with the estimate  (\ref{11}) for the problem (\ref{9}), (\ref{10}), 
since $\delta_h = \lambda_1$. Thus, we arrive at the following statement.

\begin{thm}\label{t-1}
The difference scheme (\ref{20}), (\ref{23}), (\ref{24}) is a fundamental mode exact scheme.
It is unconditionally stable for $\sigma \geq 0.5$ and the estimate (\ref{26}) holds.
\end{thm}

\section{Schemes based on Pad\'{e} approximations} 

Two-level finite difference schemes of higher approximation order 
for time-dependent linear problems can be conveniently constructed on the basis of 
Pad\'{e} approximations for the corresponding operator (matrix) exponential  function \cite{Higham2008}.  
In  the  case  of  linear  systems  of ODEs, such approximations correspond 
to various variants of the Runge-Kutta method \cite{HairerWanner2010,Butcher2008,DekkerVerwer1984}.

The solution of the problem (\ref{9}), (\ref{10}) can be represented as
\[
 w(t) = \exp(- A t) w^0 ,
 \quad 0 < t \leq T, 
\] 
and so
\[
 w(t^{n+1}) = \exp(- A \tau ) w(t^{n}) ,
 \quad n=0,1,... , N-1 . 
\]
For the corresponding difference scheme, we have
\begin{equation}\label{27}
 y^{n+1} = S y^{n} ,
 \quad n=0,1,... , N-1 ,  
\end{equation}
with the operator of transition to a new time level $S = s(- A \tau)$.
For (\ref{27}), we have the representation (\ref{15}), where now we have
\[
 y^{n+1}_k = s(\lambda_k \tau) y^{n}_k ,
 \quad k=1,2,... , M_h,  
 \quad n=0,1,... , N-1 .  
\]
It is necessary to choose one or another function $s(z)$, 
which approximates the function $\exp(-z)$ for$z \geq z_0$, where $z_0 = \lambda_1 \tau$.
To obtain an exact scheme for the fundamental mode, we must (see (\ref{16})) fulfill
$s(z_0) = \exp(-z_0)$.

The Pad\'{e} approximation of the function $\exp(-z)$ is
\begin{equation}\label{28}
  \exp(-z) = R_{lm}(z) + \mathcal{O} (z^{l+m+1}),
  \quad R_{lm}(z) = \frac{P_{lm}(z)}{Q_{lm}(z)},
\end{equation}
where $P_{lm}(z)$ and $Q_{lm}(z)$ are polynomials of degrees $l$ and $m$, respectively. 
These polynomials have \cite{BakerGraves-Morris1996} the form
\[
  P_{lm}(z) := \frac{l!}{(l+m)!} \sum_{k=0}^{l} 
  \frac{(l+m-k)!}{k! (l-k)!} (-z)^k,
\]
\[
  Q_{lm}(z) := \frac{m!}{(l+m)!} \sum_{k=0}^{m} 
  \frac{(l+m-k)!}{k! (m-k)!} z^k .
\]

At the point $z = 0$ the approximation is exact, i.e. $\exp(0) = R_{lm}(0)$.
We search the exact approximation at the point $z_1 = \eta_1 = \lambda_1 \tau$.
Taking this into account, we put
\[
   \exp(-z) \approx \exp(-z_1) R_{lm}(\widetilde{z}) ,
 \quad \widetilde{z} = z - z_1 . 
\]
The fundamental mode exact scheme based on Pad\'{e}  approximations can be
represented in the form (\ref{27}), where in view of (\ref{20}) and (\ref{28}), we have
\begin{equation}\label{29}
 S = \exp(- \lambda_1 \tau ) R_{lm}(\widetilde{A} \tau) . 
\end{equation}
In fact, this means that the transformation (\ref{25}) is realized.
Thus, we can formulate the following statement.

\begin{thm}\label{t-2}
The difference schemes (\ref{20}), (\ref{24}), (\ref{27}), (\ref{29})
based on Pad\'{e}  approximations belong to the class of fundamental mode exact schemes.
\end{thm}

Let us highlight the most interesting variants of the constructed schemes of FMES-type.
The scheme (\ref{27}), (\ref{29}) with $l = 0, m = 1$ corresponds to
the fully implicit scheme ($\sigma =1$ in(\ref{23}), (\ref{24})).
The scheme (\ref{27}), (\ref{29}) with $l = 1, m = 1$ is associated with
the Crank-Nicolson scheme ($\sigma = 0.5$ in (\ref{23}), (\ref{24})).

Only the schemes with $l = 0$, i.e.,
\[
  R_{0m}(z) = \frac{1}{Q_{0m}(z)},
  \quad P_{0m}(z) = 1
\]
are SM-stable \cite{Vabischevich2010b}. In this case, the two-level difference scheme 
for the problem (\ref{9}), (\ref{10}) is
\begin{equation}\label{30}
  \frac{\exp(\lambda_1 \tau ) y^{n+1} -y^n}{\tau} + 
  \frac{1}{\tau }(Q_{0m}(\widetilde{A} \tau) - I) y^{n+1} = 0, 
  \quad  n = 0,1, ..., N-1 ,  
\end{equation}
where the function
\[
  Q_{0m}(z) = \sum_{k=0}^{m} \frac{1}{k!} z^k 
\]
is a truncated Taylor series for $\exp(z)$. 

In the class of schemes (\ref{30}), in addition to the fully implicit scheme of 
the first-order approximation in time ($m=1$), we highlight the scheme with $m=2$.
In this case, the approximate solution is determined from
\begin{equation}\label{31}
  \frac{\exp(\lambda_1 \tau ) y^{n+1} -y^n}{\tau} + 
  \left (\widetilde{A} + \frac{\tau}{2} \widetilde{A}^2 \right ) y^{n+1} = 0, 
  \quad  n = 0,1, ..., N-1 .  
\end{equation}
It has, like the Crank-Nicolson scheme ($\sigma = 0.5$ in (\ref{23}), (\ref{24})),
the second-order approximation in time, but, as an SM-stable scheme,
it reproduces more accurately the behavior of higher modes.

\section{Numerical experiments} 

\begin{figure}[ht] 
  \begin{center}
    \begin{tikzpicture}
       \filldraw[color=blue!5] (0,0) rectangle +(5,5);
       \filldraw[color=blue!25] (0,0) rectangle +(2.5,2.5);
       \draw[color=blue!75] (0,0) rectangle +(5,5);
       \draw [->] (0,0) -- (6,0);
       \draw  (5.5,-0.25) node {$x_1$};  
       \draw [->] (0,0) -- (0,6);
       \draw  (-0.25,5.5) node {$x_2$};  
       \draw  (-0.25,-0.25) node {$0$};  
       \draw  (-0.25, 5.) node {$1$};  
       \draw  (5., -0.25) node {$1$};  
       \draw  (1.2,1.2) node {$k=10$}; 
       \draw  (3.7,3.7) node {$k=1$}; 
       \draw  (-0.5,2.5) node {$\mu=0$};  
       \draw  (5.7,2.5) node {$\mu=10$};  
       \draw  (2.5,-0.25) node {$\mu=0$};  
       \draw  (2.5,5.25) node {$\mu=10$};  
    \end{tikzpicture}
    \caption{Computational domain, coefficients of the equation and the boundary conditions}
    \label{fig:2}
  \end{center}
\end{figure}
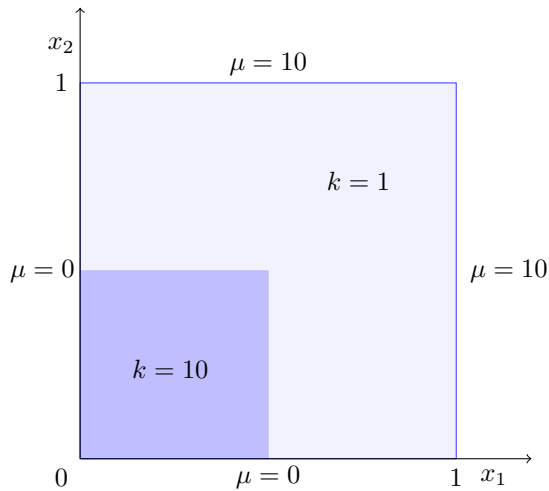

The model problem (\ref{1})--(\ref{4}) is considered in the unit square ($\Omega = [0,1]\times [0,1]$).
We assume that $c(\bm x) = c = \mathrm{const}$, and the coefficient $k(\bm x)$ in equation (\ref{1}) is discontinuous. 
Namely, in the left lower quarter of the computational domain $\Omega$,
it is equal to 10, in other parts of the region it equals 1. 
Neumann-type or Robin-type boundary conditions are specified on different parts of the boundary. 
The values of the coefficient $k(\bm x)$ and function $\mu(\bm x)$ in the boundary condition (\ref{2})
are shown in Fig.~\ref{fig:2}.
Thus, we consider the problem (\ref{3}), (\ref{4}), where
\[
 \mathcal{A} =  \bar{\mathcal{A}} + cI,
 \quad  \bar{\mathcal{A}} u :=  - \nabla (k(\bm x) \nabla u ) ,
 \quad \bm x \in \Omega . 
\]
A uniform grid in space $N_h \times N_h$ is used. Piecewise-linear finite elements on triangles are employed.

The constructed nonstandard difference schemes are based on using the known
fundamental eigenvalue $\delta = \lambda_1$.
The computation of the first eigenvalues and the corresponding eigenfunctions is
the standard problem of computational mathematics \cite{Saadbook}.
We can focus on applying well-designed algorithms and relevant open-source software.
To solve the spectral problems, we recommend the SLEPc package
(Scalable Library for Eigenvalue Problem Computations, http://slepc.upv.es/).
In this library, in particular, there are implemented the Krylov-Schur algorithm and a variant of
the Arnoldi method proposed by \cite{stewart2002krylov}. 

\begin{figure}[tbp]
  \begin{center}
    \includegraphics[width=0.85\linewidth] {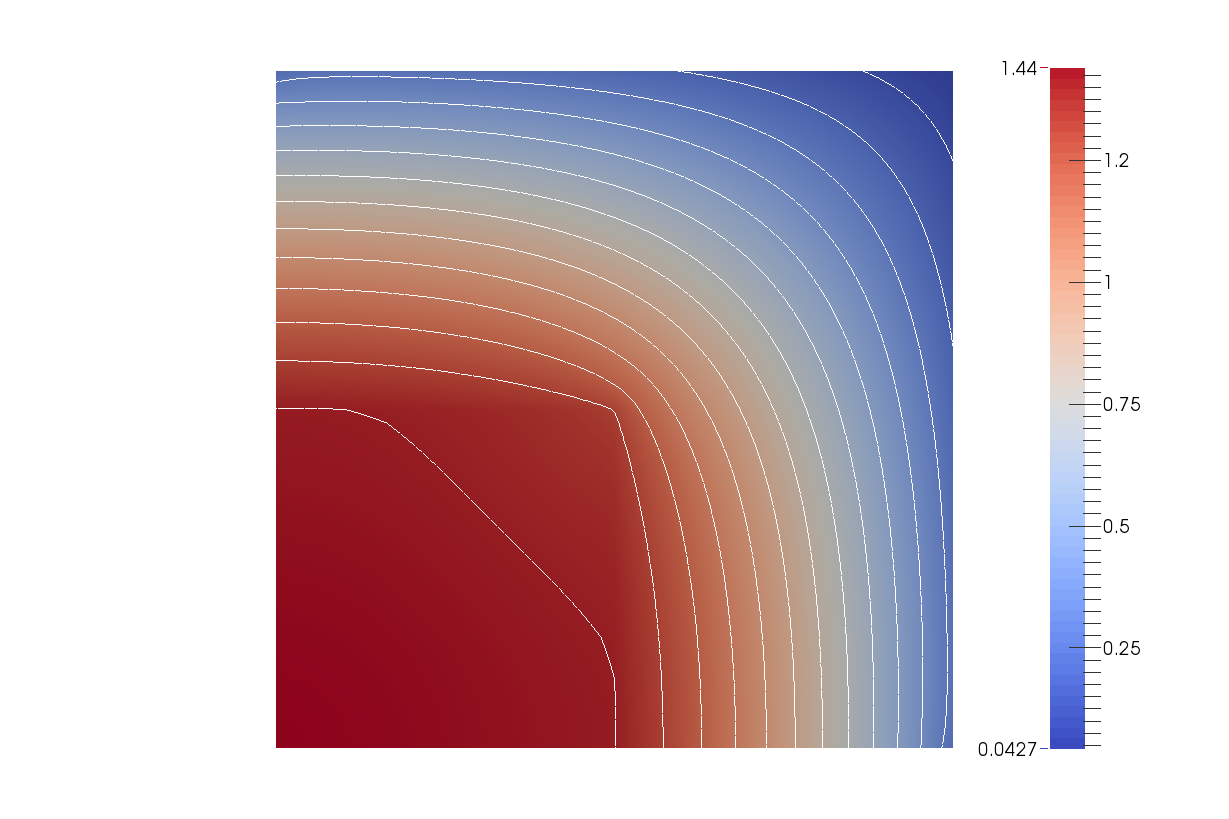}
	\caption{Fundamental eigenfunction.}
	\label{fig:3}
  \end{center}
\end{figure}

We apply, for example, the standard inverse iteration method to find the fundamental eigenvalue
of the spectral problem:
\[
 \bar{A} \varphi = \bar{\lambda} \varphi.
\]
In the simplest case, we can put $\varphi^0$ as the initial value.
At the $m+1$-st iteration, we have
\begin{equation}\label{32}
 \bar{A} \varphi^{m+1} = \varphi^m.  
\end{equation} 
For the $m+1$-st approximate values of the first eigenvalue and the corresponding eigenfunction,
we get
\[
 \bar{\lambda}_1^{m+1} = \frac{(\varphi^m,\varphi^m)}{(\varphi^{m+1},\varphi^m)} ,
 \quad  \varphi^{m+1}_1 = \frac{\varphi^{m+1}}{\|\varphi^{m+1}\|} .
\] 
For the operator $A = \bar{A} + cI$, we have the same eigenfunctions, and $\lambda_1 = \bar{\lambda}_1 + c$.
The fast convergence of the inverse iteration method for the simplest variant (\ref{32})
is presented in Table~\ref{tab-1}. The calculations are performed on the sequence of
grids derived via dividing the size of triangles by half. A small number of iterations
is enough to find the first eigenvalue and the corresponding eigenfunction.
The first eigenfunction on the grid with $N_h = 51$ is shown in Fig.~\ref{fig:3}. 
The solution value on contour levels is a multiple of $0.1$.

\begin{table}[h]
\caption{Iterative approximation to the first eigenvalue}
\label{tab-1}
\begin{center}
\begin{tabular}{rlll}
\hline
$m$ & $N_h = 26$ & $N_h = 51$ &  $N_h = 101$  \\ 
\hline
1 & 5.48146728860   & 5.41108707044   & 5.40129163974  \\
2 & 4.62435922303   & 4.54593036350   & 4.54096438823  \\
3 & 4.61225694503   & 4.53300663413   & 4.52815998124  \\
4 & 4.61203189452   & 4.53275210569   & 4.52790926702  \\
5 & 4.61202756655   & 4.53274692872   & 4.52790419871  \\
6 & 4.61202748265   & 4.53274682270   & 4.52790409555  \\
7 & 4.61202748102   & 4.53274682052   & 4.52790409345  \\
8 & 4.61202748099   & 4.53274682048   & 4.52790409341  \\
9 & 4.61202748099   & 4.53274682048   & 4.52790409340  \\
10& 4.61202748099   & 4.53274682048   & 4.52790409340  \\

\hline
\end{tabular}
\end{center}
\end{table}

\begin{figure}[tbp]
  \begin{center}
    \includegraphics[width=0.75\linewidth] {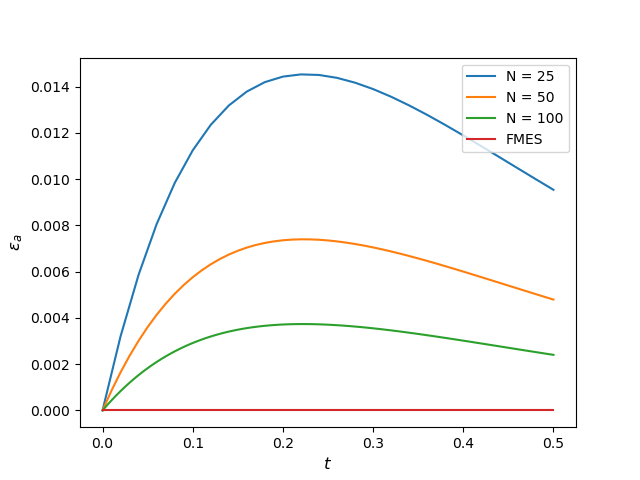}
	\caption{Time-variation the amplitude of the fundamental mode.}
	\label{fig:4}
  \end{center}
\end{figure}

The time-dependent problem is considered until $T = 0.1$.
The initial condition (\ref{4}) is taken in the form $u^0(\bm x) = 1, \ \bm x \in \Omega$.
For the basic variant, we put $c = 0$. We compare the standard fully implicit scheme 
(($\sigma =1$ in (\ref{13})) and the nonstandard scheme (\ref{23}) with the same value $\sigma$.
The accuracy of calculating the amplitude of the fundamental mode is estimated by the following value:
\[
 \varepsilon_a(t^n) = (y^n, \varphi_1) - (y^0,\varphi_1) \exp(-\lambda_1 t^n) ,
 \quad n = 0,1, ..., N .  
\] 
Numerical results for the amplitude of the fundamental mode obtained on different grids in time 
are given in Fig.~\ref{fig:4}.

\begin{figure}[tbp]
  \begin{center}
    \includegraphics[width=0.75\linewidth] {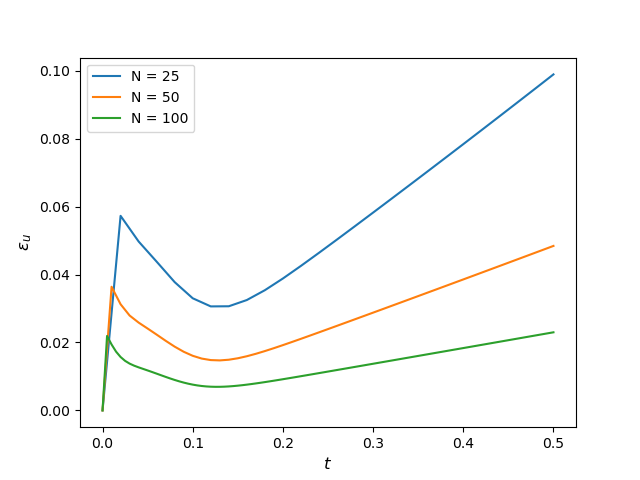}
	\caption{Accuracy of the standard fully implicit scheme.}
	\label{fig:5}
  \end{center}
\end{figure}

\begin{figure}[tbp]
  \begin{center}
    \includegraphics[width=0.75\linewidth] {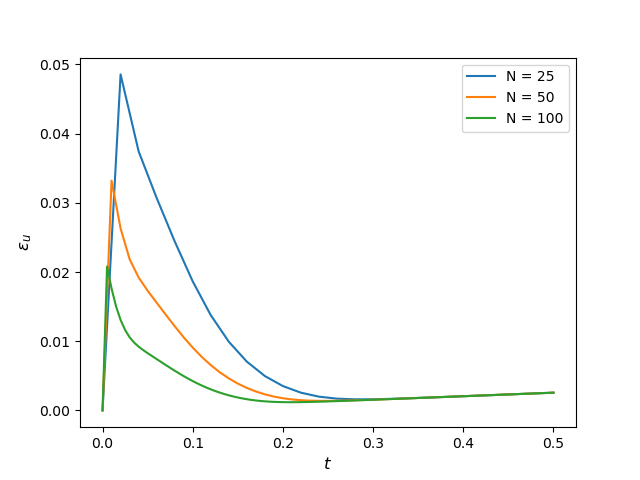}
	\caption{Accuracy of the nonstandard fully implicit scheme.}
	\label{fig:6}
  \end{center}
\end{figure}

\begin{figure}[tbp]
  \begin{center}
    \includegraphics[width=0.75\linewidth] {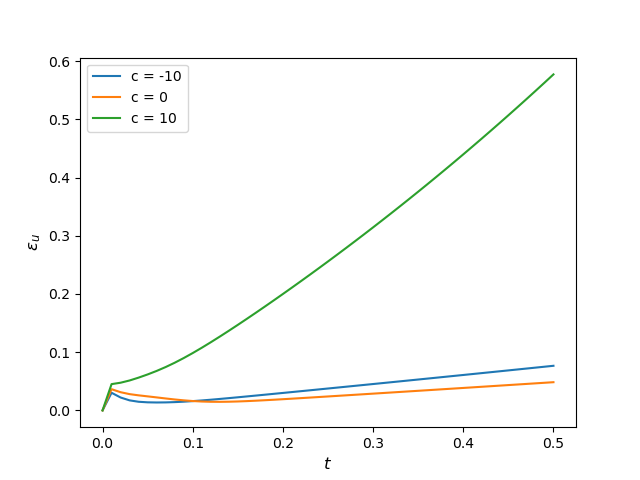}
	\caption{Accuracy of the standard fully implicit scheme for various values of $c$.}
	\label{fig:7}
  \end{center}
\end{figure}

\begin{figure}[tbp]
  \begin{center}
    \includegraphics[width=0.75\linewidth] {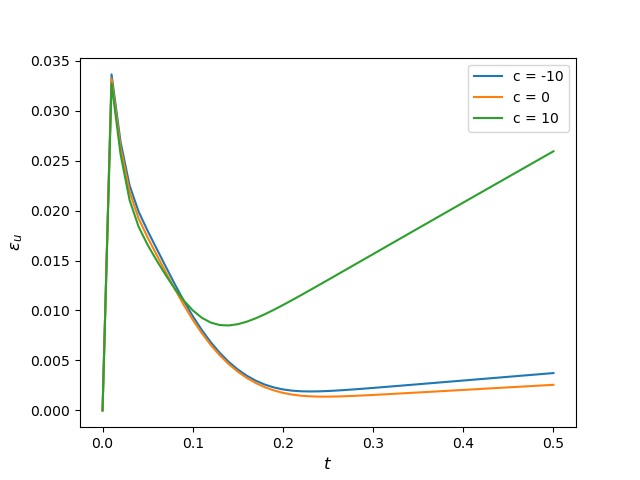}
	\caption{Accuracy of the nonstandard fully implicit scheme for various values of $c$.}
	\label{fig:8}
  \end{center}
\end{figure}

It is possible to evaluate the general error of the approximate solution without separating the fundamental mode.
As a reference solution (the benchmark solution used  for a comparison) we employ
the solution $\bar{y}$ that is obtained using the fully implicit scheme (\ref{13}), (\ref{14})
on the finest grid in time (the number of time steps is $N = 1000$).
The accuracy of the approximate solution is estimated at each time level:
\[
 \varepsilon_u(t^n) = \frac{\|y^n - \bar{y}(t^n)\|}{\|y^n\|} ,
 \quad n = 0,1, ..., N .  
\]
Numerical results for the fully implicit scheme (\ref{13}), (\ref{14}) are given in Fig.~\ref{fig:5}. 
Similar data for the nonstandard scheme (\ref{23}), (\ref{24}) for $\sigma = 1$ are shown in Fig.~\ref{fig:6}.
We observe a much higher accuracy of the approximate solution obtained using the FMES-type scheme. 
The influence of the constant $c$ is illustrated in Figs.~\ref{fig:7} and \ref{fig:8}.
Obviously, the nonstandard scheme demonstrates higher accuracy in comparison with the standard scheme.

\clearpage

\section*{Acknowledgements}

The publication was financially supported by the Ministry of Education and Science of 
the Russian Federation (the Agreement \#~02.a03.21.0008).


\end{document}